\documentclass{amsart}

\usepackage{amssymb}

\def\declaration#1#2{
\expandafter\def\csname #1\endcsname##1{\subsection{#2}\label{#1:##1}} }
\declaration{dfn}{Definition}
\declaration{xmp}{Example}
\declaration{thm}{Theorem}
\declaration{prp}{Proposition}
\declaration{crl}{Corollary}
\declaration{lem}{Lemma}
\declaration{rem}{Remark}

\def\proof{
\smallskip
\noindent{\it Proof}}

\def\ge{\geqslant}
\def\le{\leqslant}

\def\t{{\mathbb T}}

\begin{document}

\title{On Kan fibrations for Maltsev algebras}
\author{Mamuka Jibladze}
\address{Universit\'e catholique de Louvain\\
D\'epartement de math\'ematique\\
Chemin du Cyclotron 2\\
1348 Louvain-la-Neuve\\
Belgium}
\author{Teimuraz Pirashvili}
\address{A. Razmadze Mathematical Institute\\
M. Alexidze str. 1\\
Tbilisi 380093\\
Republic of Georgia}

\maketitle

It is well known that simplicial groups are Kan complexes and more generally
any surjective homomorphism of simplicial groups is  a Kan fibration. From
results in \cite{car} among other things the following statements follow:

\thm{malcevidankani}
{\sl Any simplicial model of a Maltsev theory is a Kan complex.}

\thm{kanidanmalcevi}
{\sl If $\t$ is an algebraic theory, such that
any simplicial $\t$-model is a Kan complex, then $\t$ is
a Maltsev theory.
}

\

The aim of this note is twofold, firstly to give direct proofs of these facts
without using category theory machinery and secondly to get sharper results.
In particular we prove that any surjective homomorphism of Maltsev algebras
is a Kan fibration.

A \emph{Maltsev operation} in an algebraic theory $\t$ is a ternary operation $[\_,\_,\_]$ in $\t$ satisfying the identities
$$[a,a,b]=b \textrm{ and } [a,b,b]=a.$$
An algebraic theory is called Maltsev if it possesses a Maltsev operation.
Clearly the theory of groups is a Maltsev theory by taking $[a,b,c]=ab^{-1}c.$
More generally theory of loops is a Maltsev theory; this last fact has already been used in homotopy theory -- see \cite{Kl}. Another sort
of example of Maltsev theory is the theory of Heyting algebras.

We start by proving a more strong version of Theorem \ref{thm:kanidanmalcevi}.
We let
$$
S^1= {\Delta ^1}/\partial {\Delta ^1}
$$
be the smallest simplicial model of the circle. Moreover, we let $S^1_{\t}$
be the simplicial $\t$-model, which is obtained by applying the free
$\t$-model functor degreewise to $S^1$.

\prp{strong}
{\sl  Let $\t$ be an algebraic theory such that
$S^1_{\t}$ satisfies the $(1,2)$-th Kan condition in dimension
$2$. That is, for any $1$-simplices $x_1,x_2$ with $d_1x_1=d_1x_2$ there is a
$2$-simplex $x$ with $d_1x=x_1$ and $d_2x=x_2$. Then $\t$ is a Maltsev
theory.}

\proof. Denote the unique nondegenerate 1-simplex of $S^1$ by
$\sigma$ and the unique vertex $d_0\sigma=d_1\sigma$ by $*$. So $S^1_{\t}$
in dimension zero is
the free $\t$-model generated by $*$. Similarly $S^1_{\t}$ in dimension one is the free $\t$-model
generated by $s_0*$ and $\sigma$, and in dimension two it is
the free $\t$-model
generated by $s_1s_0*$, $s_0\sigma$, and $s_1\sigma$. Since
$d_1s_0*=d_1\sigma=*$, the $(1,2)$-th Kan condition implies existence of
a 2-simplex $x$ of $S^1_{\t}$  with $d_1x=s_0*$ and $d_2x=\sigma$. This means
there is an element $x(s_1s_0*,s_0\sigma,s_1\sigma)$ in the free
$\t$-model with three generators $s_1s_0*$, $s_0\sigma$, $s_1\sigma$ such
that the equalities $x(d_1s_1s_0*,d_1s_0\sigma,d_1s_1\sigma)=s_0*$ and
$x(d_2s_1s_0*,d_2s_0\sigma,d_2s_1\sigma)=\sigma$ hold in the free
$\t$-algebra with two generators $s_0*$, $\sigma$. Applying standard
simplicial identities we see that this means $x(s_0*,\sigma,\sigma)=s_0*$
and $x(s_0*,s_0*,\sigma)=\sigma$, i.~e. that $x$ is a Maltsev operation.
\qed

\

The following theorem shows that if $\t$ is a Maltsev theory then all
surjective homomorphisms of simplicial $\t$-models are Kan fibrations which
obviously implies Theorem \ref{thm:malcevidankani}. Our proof uses exactly
the same inductive argument as the one given in \cite{cur} for simplicial
groups (see page 130 in \cite{cur}), except that we modify slightly the
initial input.

\thm{}
{\sl Any surjective homomorphism $f:X\to Y$ of simplicial models of a Maltsev
theory is a Kan fibration.
}

\proof. For $n>0$ and $0\le k\le n$, given $y\in Y_n$ with $d_iy=f(x_i)$ for
$i\ne k$, $0\le i\le n$, where $x_i$ are elements of $X_{n-1}$ with matching
faces, we have to find an $x\in X_n$ with $f(x)=y$ and $d_ix=x_i$ for $i\ne
k$. Take $w_{-1}\in f^{-1}(y)$ and then put
$$
w_j=[w_{j-1},s_jd_jw_{j-1},s_jx_j]
$$
for $0\le j<k$; moreover put $w_{n+1}=w_{k-1}$ and
$$
w_j=[w_{j+1},s_{j-1}d_jw_{j+1},s_{j-1}x_j]
$$
for $n\ge j>k$.

We then have by induction for $0\le j<k$
$$
f(w_j)=[f(w_{j-1}),s_jd_jf(w_{j-1}),s_jf(x_j)]=[y,s_jd_jy,s_jd_jy]=y
$$
and for $n\ge j>k$
$$
f(w_j)=[f(w_{j+1}),s_{j-1}d_jf(w_{j+1}),s_{j-1}f(x_j)]
=[y,s_{j-1}d_jy,s_{j-1}d_jy]=y.
$$
Furthermore for $0\le j<k$
\begin{align*}
d_iw_j&=[d_iw_{j-1},d_is_jd_jw_{j-1},d_is_jx_j]\\
&=
\begin{cases}
{[x_i,s_{j-1}d_{j-1}d_iw_{j-1},s_{j-1}d_ix_j]=}\\
{[x_i,s_{j-1}d_{j-1}x_i,s_{j-1}d_{j-1}x_i]=x_i,} & 0\le i<j,\\
{[d_jw_{j-1},d_jw_{j-1},x_j]=x_j,} & i=j
\end{cases}
\end{align*}
and then for $n\ge j>k$,
\begin{align*}
d_iw_j&=[d_iw_{j+1},d_is_{j-1}d_jw_{j+1},d_is_{j-1}x_j]\\
&=
\begin{cases}
{[x_i,s_{j-2}d_{j-1}d_iw_{j+1},s_{j-2}d_ix_j]=}&\ \\
{[x_i,s_{j-2}d_{j-1}x_i,s_{j-2}d_{j-1}x_i]=x_i,} & 0\le i<k,\\
{[x_i,s_{j-1}d_jd_iw_{j+1},s_{j-1}d_{i-1}x_j]=}&\ \\
{[x_i,s_{j-1}d_jx_i,s_{j-1}d_jx_i]=x_i,} & j<i\le n,\\
{[d_jw_{j+1},d_jw_{j+1},x_j]=x_j,} & i=j
\end{cases}
\end{align*}

Thus with $x=w_{k+1}$ one obtains $f(x)=y$ and $d_ix=x_i$ for $i\ne k$, $0\le
i\le n$, as desired. \qed

\end{document}